\documentclass[12pt,reqno,draft]{amsart}

\usepackage{amscd,amsmath}
\usepackage{amssymb}
\usepackage{amsthm}
\usepackage{enumerate}

\theoremstyle{plain}
\newtheorem{theorem}{Theorem}[section]
\newtheorem{lemma}[theorem]{Lemma}

\theoremstyle{definition}
\newtheorem{remark}[theorem]{Remark}

\newtheorem{definition}[theorem]{Definition}

\def\bu{\bullet}

\renewcommand{\epsilon}{\varepsilon}
\renewcommand{\phi}{\varphi}
\DeclareMathSymbol{\varnothing}{\mathord}{AMSb}{"3F}

\title{Bar complexes and formality of pull-backs}

\author{Steven Lillywhite}

\address{Department of Mathematics, University of Toronto,
 100 St. George St., Toronto, Ontario,
M5S 3G3} \email{sml@math.toronto.edu}


\begin{document}

\begin{abstract}
We prove a result concerning formality of the
pull-back of a fibration.  Our approach is to use bar complexes 
in the category of commutative differential graded algebras.
\end{abstract}

\maketitle

\section{Introduction}
In this note, we show that the pull-back of a fibration by a formal
map is formal.  The fibration is required to be totally non-homologous
to zero, and to be a formal map as well.  We are here referring to
the notion of formality in the setting of rational homotopy theory. 
This result extends a theorem of Vigu\'e-Poirrier, \cite{VP},
 where it is proved that the fibre of such a fibration is
a formal space.   Our proof makes use of bar complexes,
which, when 
we use a normalization due to Chen, become commutative differential
graded algebras useful for rational homotopy theory. 
We conclude with an example which generalizes a result of Baum and Smith,
\cite{BS}.
 
\section{Review of rational homotopy theory}
In this section we  briefly recall some notions from rational
homotopy theory.  References for this material are numerous and 
we mention \cite{AP}, \cite{BG}, \cite{Tanre}.

We introduce the category of commutative differential graded algebras
over a field $k$ of characteristic zero.  We assume that all algebras
are concentrated in non-negative degrees, have a differential which
raises degree by one, and are augmented.  Furthermore, we shall assume
that $H^0(A)\approx k$ for all algebras $A$.  We denote this category
by $k\mathcal{CDGA}$ and refer to objects in it as $k$CDGA's.  A morphism
of $k$CDGA's which induces an isomorphism on cohomology is called a
{\it quasi-isomorphism}.  There is a notion of homotopy between 
maps of $k$CDGA's which becomes an equivalence relation on the 
set of maps from $A_1$ to $A_2$, when the source, $A_1$, is of a special
type called a {\it KS complex}, which is basically a free algebra whose
differential respects an ordering on the generators.  

Among the KS complexes in $k\mathcal{CDGA}$
is an important class called {\it minimal} algebras, which are
essentially characterized by being free with decomposable differential.  For
every $k$CDGA $A$, there exists a minimal $k$CDGA $\mathcal M$, and
a quasi-isomorphism $\mathcal M\to A$.  Such a minimal algebra is called
a {\it minimal model} of $A$. It is unique up to isomorphism, and furthermore,
a map of algebras $f: A_1\to A_2$ determines a map $f:\mathcal M(A_1)
\to \mathcal M(A_2)$ which is unique up to homotopy.  
An algebra $A$ is called {\it formal} if there are quasi-isomorphisms of
$k$CDGA's $A\leftarrow \mathcal M(A) \to H(A)$. 
This is equivalent to demanding that there be a sequence of $k$CDGA
quasi-isomorphisms $A\leftarrow A_1\to \dots \leftarrow A_n \to H(A)$.

If $X$ is a path-connected topological space, then there is a functor
which associates to $X$ a $k$CDGA $A(X)$, known as the Sullivan-de
Rham algebra.  For a large class of spaces, including simply-connected
spaces of finite $\mathbb Q$-type (see definition below), 
the minimal model of $A(X)$
determines the rational homotopy of $X$, $\pi^\bu(X)\otimes\mathbb Q$.  We
say that a space $X$ is {\it formal} if $A(X)$ is formal.  Then for
these formal spaces, their rational
homotopy is determined by their cohomology algebras.  

\section{Bar complexes and Eilenberg-Moore theory}
In this section we shall discuss the theory of Eilenberg and Moore
concerning pull-backs of fibrations. 
For references, see \cite{Mc},  \cite{Smith}, or \cite{EMo}.
For the rest of the paper, we shall assume that all spaces are 
path-connected and of finite $k$-type (meaning that $H^n(X; k)$
is finite-dimensional for all $n\ge 1$).  We shall also assume
that all fibrations are Serre fibrations.

Let us suppose that we have a fibration $F\to E\overset{p}
{\to} B$ and a map $f: X\to
B$, so that we obtain a pull-back diagram:

\begin{equation}\label{diagram}
\begin{CD}
E_f @>{\tilde f}>> E \\
@V{\tilde p}VV   @V{p}VV \\
X @>{f}>>  B 
\end{CD}
\end{equation}

\noindent
Then the maps $f^*$ and $p^*$  make $A^\bu(X)$ and $A^\bu(E)$
(differential graded) modules over $A^\bu(B)$. Let us assume that $B$
is simply-connected.  Then a theorem of Eilenberg and Moore asserts
that there is an isomorphism 
\begin{equation}
\theta: Tor_{A^\bu(B)}(A^\bu(X), A^\bu(E))\to H^\bu(E_f).
\end{equation} 

\noindent
We may use the bar resolution to obtain a resolution of, say, $A^\bu(X)$
by $A^\bu(B)$-modules.  Since we are considering $A^\bu(-)$
to be the Sullivan-de Rham complex, we will be using Chen's
normalized bar resolution, see \cite{Ch} or \cite{GJP}.  

More specifically, the bar complex is

\begin{equation}
B(A^\bu(X), A^\bu(B), A^\bu(E))=\bigoplus_{i=0}^{\infty}
A^\bu(X)\otimes_k (sA^\bu(B))^{\otimes i}\otimes_k A^\bu(E)
\end{equation}

\noindent
where the tensor products are over the ground field $k$, and $s$ denotes
the suspension functor on graded vector spaces which lowers degree by one.  
Hence the degree of an element $(\alpha,\omega_1,\dots,\omega_k,\beta)$
is:  $deg(\alpha)+\sum_{i=1}^k (deg(\omega_i)-1)+deg(\beta)$,
where $\alpha\in A^\bu(X)$, $\omega_i\in A^\bu(B)$, and $\beta\in 
A^\bu(E)$. Actually,
the bar complex is bigraded.  We introduce the {\it bar degree}, denoted
$B(A^\bu(X), A^\bu(B), A^\bu(E))_{\bu}$.  The bar degree of an element
$(\alpha,\omega_1,\dots,\omega_k,\beta)$ is defined to be $-k$.  The
other grading is the normal tensor product grading, the degree of an
element $(\alpha,\omega_1,\dots,\omega_k,\beta)$ being 
$deg(\alpha)+\sum_{i=1}^k deg(\omega_i)+deg(\beta)$.

There are two differentials of total degree +1:

\begin{align}
d(\alpha,\omega_1,\dots,\omega_k,\beta)&=(d\alpha, \omega_1,\dots,\omega_k,
\beta)\\\notag
&+\sum_{i=1}^k(-1)^{\epsilon_{i-1}+1}
(\alpha,\omega_1,\dots,\omega_{i-1},d\omega_i,\omega_{i+1},\dots,\omega_k,\beta
)\\\notag
&+(-1)^{\epsilon_k}(\alpha,\omega_1,\dots,\omega_k,d\beta)\\
-\delta(\alpha,\omega_1,\dots,\omega_k,n)&=(-1)^{\epsilon_0}
(\alpha\omega_1,\omega_2,\dots,\omega_k,\beta)\\\notag
&+\sum_{i=1}^{k-1}(-1)^{\epsilon_i}(\alpha,\omega_1,\dots,
\omega_{i-1},\omega_i\omega_{i+1},\omega_{i+2},\dots,\omega_k,\beta)\\\notag
&+(-1)^{\epsilon_{k-1}+1}(\alpha,\omega_1,\dots,\omega_{k-1},\omega_k\beta)
\end{align}

\noindent
where 
$\epsilon_i=deg\alpha+
deg\omega_1+\dots +deg\omega_i-i$. The differential $\delta$ has
degree $+1$ with respect to the bar grading, while the differential
$d$ has degree $+1$ with respect to the tensor product grading. 
One may verify that  
$d\delta+\delta d=0$, and we put $D\overset{def}{=}
d+\delta$ to be the total differential.
With the given bigrading,
we get a double complex with the two differentials $d$ and $\delta$.
If we filter the bar complex so that we take $d$-cohomology first
in the associated spectral sequence, then we obtain the spectral 
sequence of Eilenberg and Moore. 

Chen's normalized version of this bar complex is the following.
  If $f\in A^0(B)$, let $S_i(f)$ be the operator on 
$B(A^\bu(X),A^{\bu}(B),A^\bu(E))$ defined by 
\begin{equation}
S_i(f)(\alpha,\omega_1,\dots,
\omega_k,\beta)
=(\alpha,\omega_1,\dots,\omega_{i-1},f,\omega_i,\dots,\omega_k,\beta)
\end{equation}
for
$1\le i\le k+1$.  Let $W$ be the subspace of $B(A^\bu(X),A^{\bu}(B),A^\bu(E))$
generated by the images of $S_i(f)$ and $DS_i(f)-S_i(f)D$.  Then define

\begin{equation}
\bar B(A^\bu(X),A^{\bu}(B),A^\bu(E))\overset{def}{=}
B(A^\bu(X),A^{\bu}(B),A^\bu(E))/W.
\end{equation}

\noindent  
Then $W$ is closed under $D$ and when $H^0(B)=k$ (B is connected), then
$W$ is acyclic so that $\bar B(A^\bu(X), A^\bu(B), A^\bu(E))$ 
is quasi-isomorphic
to $B(A^\bu(X), A^\bu(B), A^\bu(E))$.  Notice that in the normalized
bar complex, there are no elements of negative degree, and 
with our assumption that $B$ is simply-connected, we are
assured convergence of  the associated Eilenberg-Moore spectral sequence.
The map $\theta$ 
mentioned above is
induced by the map 
\begin{equation}
\theta: B(A^\bu(X), A^\bu(B), A^\bu(E))\to A^\bu(E_f)
\end{equation}
which sends all tensor products to zero except for 
$A^\bu(X)\otimes_k A^\bu(E)$, where the map is: $\alpha\otimes\beta\mapsto
\tilde p^*\alpha\wedge \tilde f^*\beta$.  Note that $\theta(W)=0$, so that
we get an induced map 
\begin{equation}
\theta: \bar B(A^\bu(X), A^\bu(B), A^\bu(E))\to
A^\bu(E_f).  
\end{equation}
The bar complex computes Tor, and the theorem of Eilenberg and Moore
states that this map $\theta$ is a quasi-isomorphism.

The normalized bar complex may also be augmented.  The augmentation,
$\epsilon$, maps all elements of positive total degree to zero.  The
elements of degree zero have the form $(f,g)$, where $f\in A^0(X)$
and $g\in A^0(E)$.  Then we define $\epsilon(f,g)=\epsilon_X(f)\epsilon_E(g)
=f(x_0)g(e_0)$,
where $x_0$ and $e_0$ are chosen base-points in $X$ and $E$, respectively,
and $\epsilon_X, \epsilon_E$ are the augmentations of 
$A^\bu(X), A^\bu(E)$, respectively.
If we choose base-points so that the pull-back diagram above preserves
all base-points, then $\theta$ is an augmentation preserving map.

The bar complex has a natural coalgebra structure.  Since we are
inputting $k$CDGA's to the bar complex, we also obtain a commutative
differential graded algebra structure 
on the bar complex via the shuffle product.  

More specifically, if $(a_1, \dots, a_p)$
and $(b_1, \dots, b_q)$ are two ordered sets, then a {\it shuffle}
$\sigma$ of $(a_1, \dots, a_p)$ with $(b_1, \dots, b_q)$ is a permutation
of the ordered set $(a_1, \dots, a_p, b_1, \dots, b_q)$ which preserves
the order of the $a_i$'s as well as the order of the $b_j$'s.  That is,
we demand that if $i<j$, then $\sigma(a_i)<\sigma(a_j)$ and 
$\sigma(b_i)<\sigma(b_j)$.  

We  obtain a product on $\bar B(A^\bu(X), A^\bu(B), A^\bu(E))$
by first taking the normal tensor product on the $A^\bu(X)\otimes A^\bu(E)$
factors, then taking the tensor product of this product with the 
shuffle product on the $A^\bu(B)^{\otimes i}$ factors.  As usual,
we introduce a sign $(-1)^{|\alpha||\beta|}$ whenever $\alpha$
is moved past $\beta$. The actual formula for the product is the
following.  Let $a, x\in A^\bu(X)$, $b_i, y_i\in A^\bu(B)$, and
$c, z\in A^\bu(E)$. 

\begin{equation}
(a, b_1,\dots, b_k, c)\bu(x, y_1,\dots, y_l, z)=
\sum_{\sigma}(-1)^{\eta +n_\sigma}(ax, \sigma(b_1,\dots,b_k;y_1,\dots,
y_l), cz)
\end{equation}
where $\eta=|c||x|+|x|\{|b_1|+\dots +|b_k|-k\}+|c|\{|y_1|+\dots +|y_l|-l\}$,
the sum is over all shuffles $\sigma$ of the set
$(b_1,\dots, b_k)$ with the set $(y_1, \dots, y_l)$,
and 
$$n_\sigma=\sum_{(i,j)}(|b_i|-1)(|y_j|-1)$$
where the sum is over all pairs
$(i,j)$ such that $b_i$ is moved past $y_j$ in the shuffle $\sigma$.  
One may check that this product is associative, graded commutative, and
that the differential $D$ is a derivation with respect to this product.

\begin{lemma}
The product defined above induces a product on Chen's normalized bar
complex.
\end{lemma}

\begin{proof}
We have that $\bar B(A^\bu(X), A^\bu(B), A^\bu(E))=
B(A^\bu(X), A^\bu(B), A^\bu(E))/W$.
We will show  that $W$ is an ideal.  Now $W$ is generated by
the images of $S_i(f)$ and $DS_i(f)-S_i(f)D$.  It is obvious that
$\alpha\bu S_i(f)\beta\in W$ for any $\alpha, \beta$.  Moreover,
$\alpha\bu (DS_i(f)\beta-S_i(f)D\beta)=\alpha\bu DS_i(f)\beta
-\alpha\bu S_i(f)D\beta$.  Now $\alpha\bu S_i(f)D\beta\in W$ as we
just noted.  Moreover, since $D$ is a derivation, we have
\begin{equation}\label{der}
D(\alpha\bu S_i(f)\beta)= D\alpha\bu S_i(f)\beta+ (-1)^{|\alpha|}
\alpha\bu DS_i(f)\beta
\end{equation}
Now, $\alpha\bu S_i(f)\beta\in W$, and $W$ is closed under $D$, 
so the left-hand side of \ref{der} is in $W$.  Also,
$D\alpha\bu S_i(f)\beta\in W$.  Hence, $\alpha\bu DS_i(f)\beta\in W$
as well.

\end{proof}

We have arrived at the following lemma.

\begin{lemma}\label{applemma1}
Assume that we have the pull-back diagram  \ref{diagram}, 
where $p$ is a fibration and $B$ is simply connected. Then the normalized
bar complex 
$$
\bar B(A^\bu(X), A^\bu(B), A^\bu(E))
$$
 is a $k$CDGA. 
Moreover, 
$$\theta : \bar B(A^\bu(X), A^\bu(B), A^\bu(E))\to A^\bu(E_f)$$
is a quasi-isomorphism of $k$CDGA's.  
\end{lemma}

\begin{remark}
We note that Chen's normalization is functorial.  That is, if we have
a commutative diagram of $k$CDGA's
\begin{equation}
\begin{CD}
A_2 @<<< B_2 @>>> C_2 \\
@AAA   @AAA   @AAA  \\
A_1 @<<< B_1 @>>> C_1
\end{CD}
\end{equation}
then we get a map of $k$CDGA's $\bar B(A_1, B_1, C_1)\to 
\bar B(A_2, B_2, C_2)$.
\end{remark}

The next lemma concerns quasi-isomorphisms of bar complexes.
The main idea of the proof may be found in \cite{VP}, Lemme 4.3.3, and
so we omit the proof here.

\begin{lemma}\label{barquasi}
Suppose that $A_1\leftarrow B_1\to C_1$ and $A_2\leftarrow B_2\to C_2$
are two sequences of maps of $k$CDGA's with $H^1(B_1)=0=H^1(B_2)$.
Suppose further
that $B_1$ is a KS-complex, and that we have a homotopy commutative diagram 
of the form

\begin{equation}
\begin{CD}
A_2 @<<< B_2 @>>> C_2 \\
@AAA     @AAA     @AAA  \\
A_1 @<<< B_1 @>>> C_1
\end{CD}
\end{equation}

\noindent
where the vertical arrows are all $k$CDGA quasi-isomorphisms.
Then $\bar B(A_1, B_1, C_1)$ is quasi-isomorphic to
$\bar B(A_2, B_2, C_2)$ in $k$CDGA (via a sequence of $k$CDGA
quasi-isomorphisms).

\end{lemma}

\section{Formality of pull-backs}
We can use the normalized bar complex to extend a result of Vigu\'e
-Poirrier concerning formality of the fiber of a fibration,
\cite{VP}, Th\'eor\`eme 4.4.4.  
Our proof is also shorter and more direct than in \cite{VP},
which deals with more general considerations.

\begin{definition} Suppose that $A\overset{f}{\leftarrow}B
\overset{g}{\to} C$ are morphisms of $k$CDGA's.  Then we shall
say that $f$ and $g$ are {\it compatibly formal} if there exists a
homotopy commutative diagram

\begin{equation}
\begin{CD}
A @<{f}<< B @>{g}>> C\\
@AAA  @AAA  @AAA \\
\mathcal M(A) @<<< \mathcal M(B) @>>> \mathcal M(C) \\
@VVV  @VVV  @VVV  \\
H(A) @<{f}<< H(B) @>{g}>>  H(C)
\end{CD}
\end{equation}
where the middle row are minimal models for $A, B$, and $C$, and the 
vertical arrows are quasi-isomorphisms.  We shall say that maps,
$f, g$, of spaces $X\overset{f}{\to} Y \overset{g}{\leftarrow} Z$
are {\it compatibly formal} if the the corresponding maps
$A(X)\overset{f}{\leftarrow} A(Y) \overset{g}{\to} A(Z)$ are
compatibly formal.
\end{definition}

 Consider again the pull-back diagram \ref{diagram}
where $p$ is a fibration with fiber $F$, and $B$ is simply-connected.  

\begin{theorem}\label{appprop}
Assume that $p$ and $f$ are compatibly
formal maps,
and suppose that the Serre spectral
sequence for the fibration $p$ degenerates at the $E_2$ term. Then
the pull-back $E_f$ is formal.  
\end{theorem}

\begin{proof}
By \ref{applemma1} we have a quasi-isomorphism of $k$CDGA's 
\begin{equation}
\theta : \bar B(A^\bu(X), A^\bu(B), A^\bu(E))\to A^\bu(E_f) 
\end{equation}

By the assumption of compatible formality, we have
a homotopy commutative diagram whose vertical arrows are quasi-isomorphisms
 
\begin{equation}\label{good}
\begin{CD}
A^\bu(X) @<<< A^\bu(B) @>>> A^\bu(E)\\
@AAA  @AAA  @AAA\\
\mathcal M(X) @<<< \mathcal M(B) @>>>  \mathcal M(E) \\
@VVV  @VVV  @VVV\\
H^\bu(X) @<<<  H^\bu(B) @>>> H^\bu(E)
\end{CD}
\end{equation}

\noindent
Then  we obtain a sequence of $k$CDGA quasi-isomorphisms 
\begin{equation}
\bar B(A^\bu(X), A^\bu(B), A^\bu(E))\leftarrow \dots
\to \bar B(H^\bu(X), H^\bu(B), H^\bu(E))
\end{equation}
by \ref{barquasi}.  Now the bar complex
$$\bar B(H^\bu(X), H^\bu(B), H^\bu(E))$$ \noindent
has only a single differential,
$\delta$.  Since we have assumed the Serre spectral sequence for the
fibration $p$ to degenerate at the $E_2$ term, it follows that
$H^\bu(E)$ is a free $H^\bu(B)$-module.  Thus if we grade according
to the bar degree for $\delta$,  we find that 

\begin{enumerate}
\item $H_+(\bar B(H^\bu(X), H^\bu(B), H^\bu(E))_\bu)=0$ 
\item $H_0(\bar B(H^\bu(X), H^\bu(B), H^\bu(E))_\bu)\approx 
H_\bu(\bar B(H^\bu(X), H^\bu(B), H^\bu(E))_\bu$
\end{enumerate}

\noindent
Hence the projection to cohomology
\begin{align}
\bar B(H^\bu(X), H^\bu(B), H^\bu(E))_{\bu}&\to 
\bar B(H^\bu(X), H^\bu(B), H^\bu(E))_0\notag\\
&\to H_\bu(\bar B(H^\bu(X), H^\bu(B), H^\bu(E))_\bu\notag\\
&\approx H^\bu(E_f)\notag
\end{align}
is a $k$CDGA quasi-isomorphism and $E_f$ is consequently formal.

\end{proof}

\section{An Example}
Suppose that $B$ is a simply-connected space with the property that
the cohomology of $B$ is a free $k$CDGA.  Then $B$ is a formal
space, and the cohomology, $H^\bu(B)$, is a minimal model for $A^\bu(B)$.
Let $E\overset{p}{\to} B$ be a fibration with $E$ a formal space, 
and let $X$ be a formal space with a map to $B$, $f: X\to B$.  
Then it is easy to see that $f$ and $p$ are compatibly formal.

Let $E_f$ be the pull-back of the fibration $p$ by the map $f$, as
in diagram \ref{diagram}. 
Then \ref{applemma1} and \ref{barquasi} imply that $H^\bu(E_f)$ and 
$Tor_{H^\bu(B)}(H^\bu(X), H^\bu(E))$ are isomorphic as algebras.  
This extends a result of Baum and Smith, \cite{BS}, where
this was proven by other means
for $X$ a compact globally symmetric space, and
$E=BH$, $B=BG$, for $G$ a compact, connected Lie group, and
$H\subset G$ a closed, connected subgroup.  
Moreover, if the Serre spectral sequence for the fibration
$E\overset{p}{\to} B$ degenerates at the $E_2$ term, then by
\ref{appprop}, we have that $E$ is a formal space.

\bibliographystyle{amsplain}
\bibliography{master}

\end{document}